\documentclass[a4paper,12pt]{article}

\usepackage{amssymb}
\usepackage{amsmath}
\usepackage{amstext}
\usepackage{amsthm}

\newcommand{\R}{\mathbb R}

\newcommand{\eps}{\varepsilon}

\renewcommand{\Im}{\mathrm{Im}}
\renewcommand{\Re}{\mathrm{Re}}
\newtheorem{theorem}{Theorem}
\newtheorem{lemma}[theorem]{Lemma}

\newtheorem{corollary}[theorem]{Corollary}
\theoremstyle{definition}
\newtheorem{definition}[theorem]{Definition}

\theoremstyle{remark}

\begin{document}

\title{The Nonlinear Schroedinger Equation: Existence, Stability and Dynamics
of Solitons }
\author{Vieri Benci
\thanks{Dipartimento di Matematica Applicata,
Universit\`a degli Studi di Pisa, Via F. Buonarroti 1/c, Pisa,
ITALY. e-mail: \texttt{benci@dma.unipi.it}}
\and{
Marco Ghimenti 
\thanks{Dipartimento di Matematica e Applicazioni,
Universit\`a degli Studi di Milano Bicocca,
Via Cozzi, 53, Milano, ITALY. e-mail:
\texttt{marco.ghimenti@unimib.it}}}
\and{Anna Maria Micheletti\thanks{Dipartimento di Matematica Applicata,
Universit\`a degli Studi di Pisa, Via F. Buonarroti 1/c, Pisa,
ITALY. e-mail: \texttt{a.micheletti@dma.unipi.it}}}}
\date{}
\maketitle

\begin{abstract}
In this paper we present some recent results concerning the existence, the
stability and the dynamics of solitons occurring in the nonlinear
Schroedinger equation when the parameter $h\rightarrow 0.$

We focus on the role played by the Energy and the Charge in the existence,
the stability and the dynamics of solitons. Moreover, we show that, under
suitable assumptions, the soliton approximately follows the dynamics of a
point particle, namely, the motion of its \textit{barycenter} $q(t)$
satisfies the equation
  \[ \ddot{q}(t) + \nabla V (q(t)) = H_h (t) \]
  where
  \[ \sup_{t \in \mathbb{R}} |H_h (t) | \rightarrow 0 \hspace{0.75em}
     \hspace{0.75em} \text{as} \hspace{0.75em} \hspace{0.75em} h \rightarrow
     0. \]

  {\noindent}\textbf{Mathematics subject classification}. 35Q55, 35Q51,
  37K40, 37K45, 47J35.

 {\noindent}\textbf{Keywords: }{Soliton dynamics, Nonlinear Schroedinger
  Equation, orbital stability, concentration phenomena, semiclassical limit.}
\end{abstract}
\section{Introduction}

Roughly speaking a solitary wave is a solution of a field equation whose
energy travels as a localized packet and which preserves this localization in
time.

By {\em{soliton}} we mean an {\em{orbitally stable}} solitary wave so
that it has a particle-like behavior (for the definition of orbital stability
we refer e.g. to Ref.~\cite{hylo,BBGM07,CL82,GSS87,GSS90,We85,We86}).

The aim of this paper is to review some recent results about 
the existence, the stability and the
behavior of the solitary waves relative to the equation
\begin{equation}
  \left\{ \begin{array}{l}
    ih \frac{\partial \psi}{\partial t} = - \frac{h^2}{2} \Delta \psi +
    \frac{1}{2 h^{\alpha}} W' (| \psi |) \frac{\psi}{| \psi |} + V
    (x) \psi\\
    \\
    \psi \left( 0, x \right) = \varphi_{} (x)
  \end{array} \right. \label{schr}
\end{equation}
where $\varphi_{} (x)$ is a suitable initial data.

In the first section we examine the case
$V \equiv 0$ and  $h=1$ (see Ref.~\cite{BBGM07}). 
Under suitable assumption on $W$, there exists a stationary solution of the
form $\psi \left( t, x \right) = U (x) e^{\frac{i}{h} \omega t}$, where $U$ is
a radial function decaying at infinity which solves the equation
\begin{equation}\label{staz}
  - \Delta U + W' (U) = 2 \omega U.
\end{equation}
This solution is found by a constrained minimization method that involves two
prime integrals of the motion: the Charge and the Energy. By a concentration
compactness argument it is proved that this stationary wave is
stable, so this solution is a soliton.

In the second part, we consider $V \neq 0$ and $h$ small 
(see Ref.~\cite{BGM09}). 
The stationary solution  $U (x) e^{\frac{i}{h} \omega t}$  becomes
\begin{equation}
  \psi (t, x) =  U \left( \frac{x }{h^{\beta}}
  \right) e^{
   \frac{i \omega t}{h^{\alpha + 1}}}
\end{equation}
where 

\begin{equation}\label{b>1}
\beta=1+\frac{\alpha}{2}.
\end{equation}

If $\beta > 0$,
the stationary solution concentrate as  $h \rightarrow 0$.  
Also, we can give a precise estimate of the behavior of
the Energy and the Charge.

These estimates are the key 
ingredient to study the case  $V \neq 0$. For $h$ sufficiently small, a solution
of (\ref{schr}) with initial datum
$\varphi(x)=\left[U\left(\frac{x-q_0}{h^\beta}\right)\right]
e^{\frac{1}{h} \mathbf{v} \cdot x}$ is a soliton which travels like a point
particle under the action of the potential $V$.
In fact, in the last section of this review, 
we define the barycenter $q (t)$ of a soliton (see Ref.~\cite{BGM09}) as
\begin{equation}\label{bary}
  q (t) = \frac{\displaystyle\int x| \psi (t, x) |^2  {dx}}
{\displaystyle\int | \psi (t, x) |^2 dx}
\end{equation}
and we prove that it evolves
approximatively like a point particle in a potential $V$. 
More exactly, $q(t)$ satisfies the
Cauchy problem
\begin{displaymath}
\left\{ \begin{array}{l}
     \ddot{q}_{} (t) + \nabla V (q_{} (t)) = H_h (t)\\
     q_{} (0) = q_0\\
     \dot{q}_{} (0) = \mathbf{v}
   \end{array} \right.
\end{displaymath}
where
\begin{displaymath}
\sup_{t \in \mathbb{R}} |H_h (t) | \rightarrow 0 \text{ as }h \rightarrow 0.
\end{displaymath}

In the last years there are some result about existence and dynamics of 
soliton for the Nonlinear Schroedinger Equation (see, for example 
Ref.~\cite{BL83a,BJ00,FGJS04,FGJS06,kera02,kera06,selv,squa}), 
in particular there are
results which compare the motion of the soliton with the solution of
the equation
\begin{equation}\label{equat}
\ddot{X}(t)+\nabla V(X(t))=0
\end{equation}
for $t\in (0,T]$ for some
constant $T<\infty $.

The result of Bronski and Jerrard \cite{BJ00} deals with 
a pure power nonlinearity and a bounded external
potential. The authors have shown that if the initial data is
close to $U(\frac{x-q_{0}}{h})e^{i\frac{v_{0}\cdot c}{h}}$ in a suitable
sense then the solution $\psi _{h}(t,x)$ of (\ref{schr}) satisfies for 
$t\in (0,T]$
\begin{equation}
\left\Vert \frac{1}{h^{N}}|\psi _{h}(t,x)|^{2}-
\left(\frac{1}{h^{N}}\int_{\mathbb{R}^{N}}|\psi _{h}(t,x)|^{2}dx\right)
\delta _{X(t)}\right\Vert _{C^{1\ast}}
\rightarrow 0\text{ as }h\rightarrow 0.
\end{equation}
Here $\delta _{X(t)}$ is the Dirac ``$\delta$-function'',
$C^{1\ast }$ is the dual of $C^{1}$ and $X(t)$ satisfies the equation
(\ref{equat}) with $X(0)=q_{0}$, $\dot{X}(0)=v_{0}$.

In related papers of Keraani \cite{kera02,kera06} there are slight
generalizations of the above result. Using a similar approach, Marco
Squassina \cite{squa} and Alessandro Selvitella \cite{selv} described the
soliton dynamics in an external magnetic potential.

Other results on this subject are in Ref.~\cite{FGJS04,FGJS06}. In
Ref.~\cite{FGJS04} the authors study the case of 
bounded external potential $V$.

In Ref.~\cite{FGJS06} the authors study the case of confining potential. They
assume the existence of a stable ground state solution with a null space non
degeneracy condition of the equation
\begin{equation}
-\Delta \eta _{\mu }+\mu \eta _{\mu }+W^{\prime }(\eta _{\mu })=0.
\end{equation}
The authors define a parameter $\varepsilon $ which depends on $\mu $ and on
other parameters of the problem. Under suitable assumptions they prove that
there exists $T>0$ such that, if the initial data $\psi ^{0}(x)$ is very
close to $e^{ip_{0}\cdot (x-a_{0})+i\gamma _{0}}\eta _{\mu _{0}}(x-a_{0})$
the solution $\psi (t,x)$ of problem $(P_{1})$ with initial data $\psi ^{0}$
is given by
\begin{equation}
\psi (t,x)=e^{ip(t)\cdot (x-a(t))+i\gamma (t)}\eta _{\mu (t)}(x-a(t))+w(t)
\end{equation}
with $||w||_{H^{1}}\leq \varepsilon $, $\dot{p}=-\nabla V(a)+o(\varepsilon
^{2})$, $\dot{a}=2p+o(\varepsilon ^{2})$ with $0<t<\frac{T}{\varepsilon }$
for $\varepsilon $ small.

In our paper\cite{BGM09}  we do not
require the uniqueness of the ground state solution which is, in
general, not easy to verify, and
we formulate our result such that it holds for any time $t$.

\section{Main assumptions}

In all this paper we make the following assumptions:
\begin{description}
  \item (i) the problem (\ref{schr}) has a unique solution
  \begin{equation}
    \psi \in C^0 ( \mathbb{R}, H^2 ( \mathbb{R}^N)) \cap C^1 ( \mathbb{R}, L^2
    ( \mathbb{R}^N)) \label{gv}
  \end{equation}
  (sufficient conditions can be found in Kato {\cite{Ka89}}, Cazenave
  {\cite{Ca03}}, Ginibre-Velo {\cite{GV79}}).

  \item (ii) $W : \mathbb{R}^+ \rightarrow \mathbb{R}$ is a $C^3$ function
  which satisfies:
  \begin{equation}
    W (0) = W' (0) = W^{\prime \prime} (0) = 0 \label{W}
  \end{equation}
  \begin{equation}
    |W^{\prime \prime} (s) | \leq c_1 |s|^{q - 2} + c_2 |s|^{p - 2} \text{for
    some } 2 < q \leq p < 2^{\ast}=\frac{2N}{N-2} . \label{Wp}
  \end{equation}
  \begin{equation}
    W (s) \geq - c|s|^{\nu}, \text{} c \geq 0, 2 < \nu < 2 + \frac{4}{N}
    \text{ and } s \text{ large } \label{W0}
  \end{equation}
  \begin{equation}
    \exists s_0 \in \mathbb{R}^+ \text{ such that } W (s_0) < 0 \label{W1}
  \end{equation}
  \item (iii) $V : \mathbb{R}^N \rightarrow \mathbb{R}$ is a $C^2$ function
  which satisfies the following assumptions:
  \begin{equation}
    V (x) \geq 0 ; \label{V0}
  \end{equation}
  \begin{equation}\label{Vinf}
    | \nabla V (x) | \leq V (x)^b \text{ for } |x| > R_1 > 1, b \in (0, 1) ;
  \end{equation}
  \begin{equation}
    V (x) \geq |x|^a \text{ for } |x| > R_1 > 1, a > 1. \label{Vinf1}
  \end{equation}
  \item (iv) the main assumption
\begin{equation}
 \alpha > 0
\end{equation}
\end{description}
Let us discuss the set of our assumptions:

The first assumption gives us the necessary regularity to  define the
barycenter and to prove that $q (t) \in C^2 ( \R, \R^N)$. The hypotheses on
the nonlinearity are necessary in order to have a soliton type solutions. In
particular, (\ref{Wp}) is a standard requirement to have a smooth energy
functional, (\ref{W1}) is the minimal requirement to have a {\em focusing}
nonlinearity and (\ref{W0}) is necessary to have a good minimization problem
to obtain the existence of a soliton. We require also that $V$ is a
{\em{confining}} potential (assumption (iii)). This is useful on the last
part of this paper, to prove the existence of a dynamics for the barycenter.

In our approach, the assumption $\alpha >0$
is crucial.  In fact, as we will see 
in Section
\ref{fio}, the energy $E_{h}$ of a soliton
$\psi$ is composed by two parts: the internal energy $J_{h}$ and
the dynamical energy $G.$ The internal energy is a kind of binding
energy that prevents the soliton from splitting, while the
dynamical energy is related to the motion and it is composed of
potential and kinetic energy.  We have
that (see Section \ref{hdiv1})
\begin{displaymath}
J_{h}\left( \psi\right) \cong h^{N\beta -\alpha } 
\end{displaymath}
and
\begin{displaymath}
G\left( \psi\right) \cong ||\psi||_{L^{2}}^{2}
\cong h^{N\beta }
\end{displaymath}
Then, we have that
\begin{displaymath}
\frac{G\left( \psi\right) }{J_{h}\left( \psi \right) }\cong
h^{\alpha  }
\end{displaymath}
So the assumption $\alpha  >0$ implies that, for $h\ll 1,$
$G\left(\psi \right) \ll J_{h}\left( \psi\right)$, namely
the internal energy is bigger than the dynamical energy. This is
the fact that guarantees the existence and the stability of the
travelling soliton for any time.

\subsection{Notations}

In the next we will use the following notations:
\begin{eqnarray*}
  &  & \Re (z), \Im (z) \text{are the real and the imaginary part of } z\\
  &  & B (x_0, \rho) =\{x \in \mathbb{R}^N : |x - x_0 | \leq \rho\}\\
  &  & S_{\sigma} =\{u \in H^1 : ||u||_{L^2} = \sigma\}\\
  &  & J_h^c =\{u \in H^1 : J_h (u) < c\}\\
  &  & \partial_t \psi = \frac{\partial}{\partial t} \psi\\
  &  & m=m_{\sigma^2} := \inf_{u \in H^1, \int u^2 = \sigma^2} J (u)\\
  &  & \beta=1+\frac{\alpha}{2}
\end{eqnarray*}

\section{General features of NSE}

Equation (\ref{schr}) is the Euler-Lagrange equation relative to the
Lagrangian density
\begin{equation}
  \mathcal{L} = \Re(ih \partial_t \psi \overline{\psi}) - \frac{h^2}{2} \left|
  \nabla \psi \right|^2 - W_h (\psi) - V (x) \left| \psi \right|^2
  \label{lagr}
\end{equation}
where, in order to simplify the notation we have set
\begin{equation}
W_h (\psi) = \frac{1}{h^{\alpha}} W (\left| \psi
   \right|) 
\end{equation}
Sometimes it is useful to write $\psi$ in polar form
\begin{equation}
  \psi (t, x) = u (t, x) e^{iS (t, x) / h} . \label{polar1}
\end{equation}
Thus the state of the system $\psi$ is uniquely defined by the couple of
variables $(u, S)$. Using these variables, the action $\mathcal{S =} \int
\mathcal{L} dxdt$ takes the form
\begin{equation}
  \mathcal{S} (u, S) = - \int \left[ \frac{h^2}{2} \left| \nabla u \right|^2 +
  W_h (u) + \left( \partial_t S + \frac{1}{2} \left| \nabla S \right|^2 + V
  (x) \right) u^2 \right] dxdt \label{polarAZ}
\end{equation}
and equation (\ref{schr}) becomes:
\begin{equation}
  - \frac{h^2}{2} \Delta u + W_h' (u) + \left( \partial_t S + \frac{1}{2}
  \left| \nabla S \right|^2 + V (x) \right) u = 0 \label{Sh1}
\end{equation}
\begin{equation}
  \partial_t \left( u^2 \right) + \nabla \cdot \left( u^2 \nabla S \right) = 0
  \label{Sh2}
\end{equation}

\subsection{The first integrals of NSE\label{fio}}

Noether's theorem states that any invariance under a one-parameter group of the
Lagrangian implies the existence of an integral of motion (see e.g.
{Gelfand-Fomin\cite{Gelfand}}).

Now we describe the first integrals which will be relevant for this paper,
namely the energy and the ``hylenic charge''.
\begin{description}
  \item[Energy] The energy, by definition, is the quantity which is preserved
  by the time invariance of the Lagrangian; it has the following form
  \begin{equation}
    E_h (\psi) = \int \left[ \frac{h^2}{2} \left| \nabla \psi \right|^2 + W_h
    (\psi) + V (x) \left| \psi \right|^2 \right] dx.
  \end{equation}
  Using (\ref{polar1}) we get:
  \begin{equation}
    E_h (u,S) = \int \left( \frac{h^2}{2} \left| \nabla u \right|^2 + W_h (u)
    \right) dx + \int \left( \frac{1}{2} \left| \nabla S \right|^2 + V (x)
    \right) u^2 dx \label{Schenergy}
  \end{equation}
  Thus the energy has two components: the \textit{internal energy} (which,
  sometimes, is also called \textit{binding energy})
  \begin{equation}
    J_h (u) = \int \left( \frac{h^2}{2} \left| \nabla u \right|^2 + W_h (u)
    \right) dx \label{j}
  \end{equation}
  and the \textit{dynamical energy}
  \begin{equation}
    G (u, S) = \int \left( \frac{1}{2} \left| \nabla S \right|^2 + V (x)
    \right) u^2 dx \label{g}
  \end{equation}
  which is composed by the \textit{kinetic energy} $\frac{1}{2} \int \left|
  \nabla S \right|^2 u^2 dx$ and the \textit{potential energy} $\int V (x)
  u^2 dx$.
  \item[Hylenic charge] Following Ref. \cite{hylo} 
the {\em{hylenic charge}},
  is defined as the quantity which is preserved by by the invariance of the
  Lagrangian with respect to the action
  \[ \psi \mapsto e^{i \theta} \psi . \]
  For equation (\ref{schr}) the charge is nothing else but the $L^2$ norm,
  namely:
  \[ \mathcal{H} (\psi) = \int \left| \psi \right|^2 dx = \int u^2 dx. \]
  \item[Momentum] If $V = 0$ the Lagrangian is also invariant by translation.
  In this case we have the conservation of the {\em momentum}
  \begin{equation}\label{momentum}
   P_j (\psi) = h \Im \int \psi_{x_j} \bar{\psi} {dx}, \ \ j = 1, 2, 3
  \end{equation}
  hence we have the first
  Newton law for the barycenter.
\end{description}

\section{The case $h = 1, V = 0$}\label{hdiv1}

In this section we present some results contained in Ref. \cite{BBGM07}.  
We minimize the internal energy $J (u)$ on the constraint $\{u
\in H^1 : \|u\|_{L^2} = \sigma\}$ for some $\sigma$ fixed. 
If $U$ is the minimizer 
and if $2\omega$ is the Lagrange multiplier associated to $U$, 
$\psi(t,x)=U(x)e^{i\omega t}$ is
a stationary solution of (\ref{schr}).

We get the following result

\begin{lemma}
\label{mainradiale} { Let $W$ satisfy (\ref{Wp}), (\ref{W0}) and
(\ref{W1}). Then, $\exists\ \bar \sigma$ such that $\forall\
\sigma>\bar\sigma$ there exists $ \bar u\in H^1$ satisfying
\begin{equation*}
J(\bar u)=m_{\sigma^2}:= \inf_{\{v\in H^1,\ ||v||_{L^2}=\sigma\}}J(v),
\end{equation*}
with $||\bar u||_{L^2}=\sigma$. Then, there exist $\omega$ and $\bar
u $ that solve (\ref{staz}), with $ \omega<0$ and $\bar u $ positive
radially symmetric.}
\end{lemma}

In order to have stronger results, we can replace (\ref{W1}) with
the following hypothesis
\begin{equation}  \label{W2}
W(s)<-s^{2+\epsilon},\ 0<\epsilon<\frac 4N\text{ for small }s.
\end{equation}
In this case we find the following results concerning the
existence of the minimizer of $J(u)$ for any $\sigma$.

\begin{corollary}
\label{cor1} { If (\ref{Wp}), (\ref{W0}) and (\ref{W2}) hold,
then for all $ \sigma$, there exists $\bar u\in H^1$, with 
$||\bar u||_{L^2}=\sigma$, such that
\begin{equation*}
J(\bar u)= \inf_{\{v\in H^1, ||v||_{L^2}=\sigma\}}J(v).
\end{equation*}}
\end{corollary}
In particular, for $N=3$ we have

\begin{corollary}
\label{cor2} { Let $N=3$. If (\ref{Wp}) and (\ref{W0}) hold and
$W\in C^{3}$, with $W'''(0)<0$, then for all
$\sigma $, there exists $\bar{u} \in H^{1}$ with
$||\bar{u}||_{L^{2}}=\sigma $ such that
\begin{equation*}
J(\bar{u})=\inf_{\{v\in H^{1},||v||_{L^{2}}=\sigma\}}J(v).
\end{equation*}}
\end{corollary}

We sketch briefly the steps of the proof for Lemma \ref{mainradiale}.
\begin{itemize}
\item[{\em Step 1:}] If $W$ satisfies (\ref{W1}) then
 $m_{\sigma^2}$:=$\inf_{S_{\sigma}} J (u)<0$
\item[{\em Step 2:}] If $W$ satisfies (\ref{W0}) then
 $m_{\sigma^2}> - \infty$, any
  minimizing Palais Smale sequence $u_n$ is bounded in $H^1$ and the 
  Lagrange multipliers $\omega_n$ associated to $u_n$ are bounded in $\R$.
\item[{\em Step 3:}] Any minimizing Palais Smale sequence converges in $H^1$ 
 to a minimizer.
\end{itemize}
We point  out that (\ref{W1}) is a fundamental requirement for the existence
of a minimizer. In fact, if $W \geq 0$, then by Pohozaev identity we can prove
that $U \equiv 0$ is the unique radial solution of (\ref{staz}). 

Concerning the stability of stationary solution we set
\begin{equation}
  S =\{U (x) e^{i \theta}, \theta \in S^1, \|U\|_{L^2} = \sigma, J (U) =
  m_{\sigma^2} \}
\end{equation}
\begin{definition}
  $S$ is orbitally stable if
\begin{displaymath}
 \forall \eps,\  \exists \delta > 0 \text{ s.t. } \forall 
\varphi \in H^1 ( \R^N),\
     \inf_{u \in S} \| \ |\psi_0| - u\|_{H^1} < \delta \Rightarrow
\end{displaymath}
\begin{displaymath}
 \forall t \inf_{u \in S} \| \ |\psi (t, \cdot)| - u\|_{H^1} < \eps
\end{displaymath}
  where $\psi$ is the solution of (\ref{schr}) with initial data $\varphi$.
\end{definition}
Using concentration compactness\cite{Li84a,Li84b} arguments we prove the 
following (see Ref.~\cite{BBGM07}, Sect. 3)
\begin{theorem}
  Let $W$ satisfy (\ref{Wp}), (\ref{W0}) and (\ref{W1}). Then $S$ is orbitally
  stable.
\end{theorem}

This variational approach can be successfully used to find stable solitary
waves for the nonlinear Klein Gordon equation
\begin{equation}
  \Box \psi = W' (| \psi |) \frac{\psi}{| \psi |} .
\end{equation}
Again, the crucial assumption to obtain solitons is (\ref{W1}) 
(see Ref.~\cite{BBBM}  for details).

We obtain a concentration result for a
minimizer $U$ crucial for this work (see Ref.~\cite{BGM09}).

\begin{lemma}
  \label{conc1}For any $\varepsilon > 0$, there exists an $\hat{R} = \hat{R}
  (\varepsilon)$ and a $\delta = \delta (\varepsilon)$ such that, for any $u
  \in J^{m + \delta} \cap S_{\sigma}$, we can find a point $\hat{q} = \hat{q}
  (u) \in \mathbb{R}^N$ such that
  \begin{equation}
    \frac{1}{\sigma^2} \int _{\mathbb{R}^N \smallsetminus B ( \hat{q},
    \hat{R})} u^2 (x) dx < \varepsilon .
  \end{equation}
\end{lemma}
We give a sketch of the proof.
\begin{proof}Firstly we prove that for any $\varepsilon > 0$, there
exists a $\delta$ such that, for any $u \in J^{m + \delta} \cap S_{\sigma}$,
we can find a point $\hat{q} = \hat{q} (u) \in \mathbb{R}^N$ and a radial
ground state solution $U$ such that
\begin{equation}
  ||u (x) - U (x - \hat{q}) ||_{H^1} \leq \varepsilon . \label{31}
\end{equation}
We argue by contradiction: if (\ref{31}) do not hold, we can construct a
minimizing sequence which not converge. At this point, given $\varepsilon$,
there exist a point $\hat{q} = \hat{q} (u) \in \mathbb{R}^N$ and a radial
ground state solution $U$ such that
\begin{equation}
  \label{formconc3} u (x) = U (x - \hat{q}) + w \text{ and } ||w||_{H^1} \leq
  C \varepsilon .
\end{equation}
Now, we choose $\hat{R}$ such that
\begin{equation}
  \label{eqconc2} \frac{1}{\sigma^2} \int _{\mathbb{R}^N \smallsetminus B (0,
  \hat{R})} U^2 (x) dx < C \varepsilon
\end{equation}
for all $U$ radial ground state solutions. This is possible 
because if $U$ is a minimizer for $J$ constrained on $S_\sigma$, 
then there exists two constants $C,R$, not depending on $U$ such that 
\begin{equation*}
U(x)<Ce^{-|x|}\text{ for }|x|>>R.
\end{equation*}
By this fact we get the claim.
\end{proof}

We remark that, depending on the nonlinearity $W$, it is
possible that the minimizer of the constrained problem is not unique.
Anyway, by Lemma
\ref{conc1}, $\hat{R}$ does not depend on the minimizer.

\section{The case $h$ small enough}

We present now the main results contained in Ref.~\cite{BGM09}.
We recall some inequalities which are useful in the following. Let it be 
\begin{equation*}
u (x) : = v \left( \frac{x}{h^{\beta}} \right).
\end{equation*}
We have
\begin{equation*}
||u||_{L^2}^2 =  \int v \left( \frac{x}{h^{\beta}} \right)^2
   dx = h^{N \beta } \int v \left( \xi \right)^2 d \xi = 
h^{N \beta} ||v||_{L^2}^2 .
\end{equation*}
and
\begin{equation}\label{Jresc}
\begin{split}
  J_h (u) & =  \int \frac{h^2}{2} | \nabla u|^2 + \frac{1}{h^{\alpha}} 
  W ( u) dx =\\
  & =  \int \frac{h^{2 }}{2} \left| \nabla_x v \left(
  \frac{x}{h^{\beta}} \right) \right|^2 + \frac{1}{h^{\alpha}} 
  W \left( v \left( \frac{x}{h^{\beta}} \right) \right) dx =\\
  & =  \int \frac{h^{N \beta + 2  - 2 \beta}}{2} \left|
  \nabla_{\xi} v \left( \xi \right) \right|^2 + h^{N \beta - \alpha }
  W \left( v \left( \xi \right) \right) d \xi =\\
  & =  h^{N \beta - \alpha } \int \frac{1}{2} \left| \nabla_{\xi} v
  \left( \xi \right) \right|^2 + W \left( v \left( \xi \right) \right) d \xi =
  h^{N \beta - \alpha } J_1 (v) .
\end{split}
\end{equation}

We give now some results about the concentration property of the solutions
$\psi (t,x)$ of the problem (\ref{schr}). Given $K>0$, $h>0$, we
put
\begin{equation}
B_{h}^{K}=\left\{
\begin{array}{c}
\varphi(x)=\psi (0,x)=u_{h}(0,x)e^{\frac{i}{h}S_{h}(0,x)}\text{ } \\
\text{with }u_{h}(0,x)=\left[ (U+w)\left( \frac{x}{h^{\beta }}
\right) \right] \\
\\
U\text{ is a minimizer of }J\text{ constrained on }S\sigma \\
\\ 
||U+w||_{L^{2}}=||U||_{L^{2}}=\sigma \text{ and }
J(U+w)\leq m+Kh^{\alpha} \\
\\
||\nabla S_{h}(0,x)||_{L^{\infty }}\leq K\text{ for all }h \\
\\
\int_{\mathbb{R}^{N}}V(x)u_{h}^{2}(0,x)dx\leq Kh^{N\beta -2\alpha }
\end{array}
\right\}.  \label{bkqh}
\end{equation}
Considering the set  $B_{h}^{K}$ as the admissible initial data set, we get
\begin{theorem}
\label{teoconcinf} Assume $V\in L_{\text{loc}}^{\infty }$ and (\ref{V0}).
Fix $K>0$, $q\in \mathbb{R}^{N}$. Let $\alpha >0 $.

For all $\varepsilon >0$, there exists $\hat{R}>0$ and $h_{0}>0$ such that,
for any $\psi (t,x)$ solution of (\ref{schr}) with initial data $\psi
(0,x)\in B_{h}^{K}$ with $h<h_{0}$, and for any $t$, there exists
$\hat{q}_{h}(t)\in \mathbb{R}^{N}$ for which
\begin{equation}
\frac{1}{||\psi (t,x)||_{L^{2}}^{2}}
\int\limits_{\mathbb{R}^{N}\smallsetminus B(\hat{q}_{h}(t),\hat{R}h^{\beta })}
|\psi(t,x)|^{2}dx<\varepsilon .
\end{equation}
Here $\hat{q}_{h}(t)$ depends on $\psi (t,x)$.
\end{theorem}

We give the proof because it is simple and quite interesting.
\begin{proof} By the conservation law, the energy $E_{h}(\psi
(t,x))$ is constant with respect to $t$. Then we have
\begin{eqnarray*}
E_{h}(\psi (t,x))&=&E_{h}(\psi (0,x))\\
&=&J_{h}(u_{h}(0,x))+\int_{\R^N}u_{h}^2(0,x)
\left[\frac{|\nabla S_{h}(0,x)|^2}2+V(x)\right]dx\\
&\leq&J_{h}(u_{h}(0,x))+\frac K2\sigma^2h^{N\beta}+Kh^{N\beta}\\
&=&h^{N\beta-\alpha}J\left(U+w\right)+C h^{N\beta}
\end{eqnarray*}
where $C$ is a suitable constant. Now, by rescaling, and using that
$\psi (0,x)\in B^{K,q}_h$, we obtain
\begin{eqnarray}
\nonumber E_{h}(\psi (t,x))
&\leq&h^{N\beta-\alpha}J(U+w)+Ch^{N\beta}\\
&\leq&
h^{N\beta-\alpha}(m+Kh^{\alpha})+Ch^{N\beta}
\label{eq77bis}\\
\nonumber
&=&h^{N\beta-\alpha}(m+Kh^{\alpha}+Ch^{\alpha})=
h^{N\beta-\alpha}\big(m+h^{\alpha}C_1\big)
\end{eqnarray}
where $C_1$ is a suitable constant. Thus
\begin{eqnarray}
\nonumber
J_{h}(u_{h_n}(t,x))&=&E_h(\psi (t,x))-G(\psi (t,x)\\
&=&E_h(\psi (t,x))-\int_{\R^N} \left[\frac{|\nabla
S_h(t,x)|^2}2+V(x)\right]u_h(t,x)^2dx\nonumber
\\
\label{Jh} &\leq&
h^{N\beta-\alpha}\big(m+h^{\alpha}C_1\big)
\end{eqnarray}
because $V\geq0$. By rescaling the inequality (\ref{Jh}) we get
\begin{equation}
J\big(u_h(t,h^\beta x)\big)\leq m+h^{\alpha}C_1
\end{equation}
So, if $\alpha>0$, for $h$ small by a simple argument and 
Lemma \ref{conc1} we get the claim.
\end{proof}
Roughly speaking we have that $J_h(\psi)\cong h^{N\beta-\alpha}$ 
and  $G(\psi)\cong h^{N\beta}$ and this is the key of the proof.

To simplify in the following we take an initial data of the type
\begin{equation}\label{datoiniziale}
\varphi(x)=U\left(\frac{x-q}{h^\beta}\right)e^{i\mathbf{v}\cdot x},
\end{equation}
where $q$, $\mathbf{v}$ are fixed. Obviously $\varphi(x)\in B_h^K$
for some $K$.

\subsection{Existence and dynamics of barycenter}

We recall the definition of barycenter of $\psi$
\begin{equation}
  q_h (t) = \frac{\displaystyle\int_{\mathbb{R}^N} x| \psi (t, x) |^2dx}
  {\displaystyle\int_{\mathbb{R}^N} | \psi (t, x) |^2 dx} . \label{eqbar}
\end{equation}
The barycenter is not well defined for all the functions $\psi \in H^1 (
\mathbb{R}^N)$. Thus we need the following result:

\begin{theorem}
  \label{din}Let $\psi (t, x)$ be a global solution of (\ref{schr}) such that
  $\psi (t, x) \in C ( \mathbb{R}, H^2 ( \mathbb{R}^N)) \cap C^1 ( \mathbb{R},
  L^2 ( \mathbb{R}^N))$ with initial data $\psi (0, x)$ such that
  \[ \int_{\mathbb{R}^N} |x|| \psi (0, x) |^2 dx < + \infty . \]
  Then the map $q_h (t) : \mathbb{R} \rightarrow \mathbb{R}^N,$ given by
  (\ref{eqbar}) is $C^2 ( \R, \R^N)$ and it holds
  \begin{equation}
    \dot{q}_h (t) = \frac{\Im \left( h \int_{\mathbb{R}^N} \bar{\psi} (t, x)
    \nabla \psi (t, x) dx \right)}{|| \psi (t, x) ||_{L^2}^2} . \label{pina}
  \end{equation}
  \begin{equation}
    \ddot{q}_h (t) = \frac{\int_{\mathbb{R}^N} V (x) \nabla | \psi (t, x) |^2
    dx}{|| \psi (t, x) ||_{L^2}^2} . \label{pino}
  \end{equation}
\end{theorem}

We have the following corollary

\begin{corollary}
\label{dincor}Assume (\ref{Vinf}) and the assumptions of the previous
theorem; then
\begin{equation}
\ddot{q}_{h}(t)=-\frac{\displaystyle\int_{\mathbb{R}^{N}}\nabla V(x)|\psi
(t,x)|^{2}dx}{||\psi (t,x)||_{L^{2}}^{2}}.  \label{pinu}
\end{equation}
\end{corollary}

\section{The final result}

\subsection{Barycenter and concentration point}

We have two quantities which describe the properties of the travelling soliton:
the concentration point $\hat q$ and the barycenter $q$. If we want to describe
the particle-like behavior of the soliton the concentration point $\hat q$ seems to be
the natural indicator: it localize at any time $t$ the center of a ball which
contains the larger part of the soliton.
Unfortunately we do not have any control on the smoothness of 
$\hat q(t)$ (indeed $\hat q$
is nor uniquely defined). The barycenter, at the contrary, 
for a very large class of solutions has
the required regularity, and the equation (\ref{pinu}) is very similar 
to the equation of the motion
we want to obtain. In this paragraph, we estimate the distance 
between the concentration point
and the barycenter of a solution $\psi (t, x)$ for a potential satisfying
hypothesis (\ref{V0}) and (\ref{Vinf1}), say a {\em confining} potential.

The assumption (\ref{Vinf1}) is necessary if we want to identify
the position of the soliton with the barycenter. Let
us see why. Consider a soliton $\psi (x)$ and a perturbation
\begin{equation*}
\psi _{d}(x)=\psi (x)+\varphi \left( x-d\right) ,\ d\in \mathbb{R}^{N}
\end{equation*}
Even if $\varphi \left( x\right) \ll \psi (x),$ when $d$ is very large, the
``position'' of $\psi (x)$ and the barycenter of $\psi _{d}(x)$ are far from
each other. In Lemma \ref{lemmabar5}, we shall prove that this situation
cannot occur provided that (\ref{Vinf1}) hold. In a paper \cite{BGMwp} 
in preparation, we
give a more involved notion of barycenter of the soliton and we will be able
to consider other situations.

Hereafter, fixed $K > 0$, we assume that $\psi (t, x)$ is a global solution of
the Schroedinger equation (\ref{schr}), $\psi (t, x) \in C ( \mathbb{R}, H^1)
\cap C^1 ( \mathbb{R}, H^{- 1})$, with initial data $\psi (0, x)
\in B_h^{K}$ with $B_h^{K}$ given by (\ref{bkqh}). We start with some technical lemma.

\begin{lemma}
  \label{lemmabar1}There exists a constant $L > 0$ such that
  \[ 0 \leq \frac{1}{h^{N \beta - 2 \alpha}} \int_{\mathbb{R}^N} V (x) u_h^2
     (t, x) dx \leq L \ \ \forall t \in \mathbb{R} . \]
\end{lemma}
\noindent The proof follows by estimating the energy.

\begin{lemma}
  \label{lemmabar3}There exists a constant $K_1$ such that
  \[ |q_h (t) | \leq K_1 \text{ for } t \in \mathbb{R} . \]
\end{lemma}

\noindent The proof follows by Lemma \ref{lemmabar1} and by (\ref{Vinf1}). Furthermore,
we can choose $R_2$ such that
\begin{equation*}
\frac{\displaystyle\int_{|x|\geq R_{2}}u_{h}^{2}(t,x)dx}
{\displaystyle\int_{\mathbb{R}^{N}}u_{h}^{2}(t,x)dx}\leq
\frac12
\end{equation*}

\begin{lemma}
  \label{lemmabar5}Given $0 < \varepsilon < 1 / 2$, and $R_2$ as in the
  previous lemma.

  We get
  \begin{enumerate}
    \item $\sup_{t \in \mathbb{R}} | \hat{q}_h (t) | < R_2 + \hat{R}
    (\varepsilon) h^{\beta} < R_2 + 1$, for all $h < \bar{h}$ and $\delta <
    \bar{\delta}$ small enough.

    \item $\sup_{t \in \mathbb{R}} \big| q_h (t) - \hat{q}_h (t) \big| <
    \frac{3 L}{\sigma^2 R_3^{a - 1}} + 3 R_3 \varepsilon + \hat{R}
    (\varepsilon) h^{\beta}$, for any $R_3 \geq R_2$, and for all $h$ small
    enough.
  \end{enumerate}
\end{lemma}

\noindent The hardest part of the proof is the estimate of
\begin{equation*}
  I_1 = \frac{\left| \int_{\mathbb{R}^N \smallsetminus B (0, R_3)}
  \left( x - \hat{q}_h (t) \right) u_h^2 (t, x) dx
  \right|}{\int_{\mathbb{R}^N} u_h^2 (t, x) dx}.
\end{equation*}
Using (\ref{Vinf1}) and the previous estimates we can conclude.

We notice that $R_1, R_2$ and $R_3$ defined in this section do not depend on
$\varepsilon$.

\subsection{Equation of the travelling soliton}

We prove that the barycenter dynamics is approximatively that of a point
particle moving under the effect of an external potential $V (x)$.

\begin{theorem}
  \label{mainteoloc}Assume (i)-(iv). Given $K > 0$, let
  $\psi (t, x) \in C (\mathbb{R}, H^2) \cap C^1 ( \mathbb{R}, H^1)$
  be a global solution of
  equation (\ref{schr}), with initial data in $B_h^{K}$, $h < h_0$. Then we
  have
  \begin{equation}
    \ddot{q_h} (t) + \nabla V (q_h (t)) = H_h (t)
  \end{equation}
  with $||H_h (t) ||_{L^{\infty}}$ goes to zero when $h$ goes to zero.
\end{theorem}

\textbf{Proof.} We know by Theorem \ref{din}, that
\begin{equation}
  \ddot{q_h} (t) +
  \frac{\displaystyle\int_{\mathbb{R}^N} \nabla V (x) u_h^2 (t, x) dx}
  {\displaystyle\int_{\mathbb{R}^N} u_h^2 (t, x) dx} = 0 \label{eqtrav1}
\end{equation}
Hence we have to estimate the function
\begin{equation}
  H_h (t) = [\nabla V ( \hat{q}_h (t)) - \nabla V (q_h (t))] +
  \frac{\int_{\mathbb{R}^N} [\nabla V (x) - \nabla V ( \hat{q}_h (t))] u_h^2
  (t, x) dx}{\int_{\mathbb{R}^N} u_h^2 (t, x) dx} . \label{eqtrav2}
\end{equation}

By Lemma \ref{lemmabar3} and Lemma \ref{lemmabar5} we get
\begin{eqnarray*}
  \big| \nabla V ( \hat{q}_h (t)) - \nabla V (q_h (t)) \big| & \leq &
  \max_{{{\footnotesize \begin{array}{l}
    i, j = 1, \ldots, N\\
    | \tau | \leq K_1 + R_2 + 1
  \end{array}}}} \left| \frac{\partial^2 V (\tau)}{\partial x_i \partial x_j}
  \right| | \hat{q}_h (t) - q_h (t) | \leq \\
  & \leq & M \left[ \frac{3 L}{\sigma^2 R_3^{a - 1}} + 3 R_3 \varepsilon +
  \hat{R} (\varepsilon) h^{\beta} \right],  \label{formH1}
\end{eqnarray*}
for any $R_3 \geq R_2$ and some $M>0$.

To estimate
\[ \frac{\int_{\mathbb{R}^N} [\nabla V (x) - \nabla V ( \hat{q}_h (t))] u_h^2
   (t, x) dx}{\int_{\mathbb{R}^N} u_h^2 (t, x) dx} \]
we split the integral three parts.
\begin{eqnarray*}
  L_1 & = & \frac{\int_{B ( \hat{q}_h (t), \hat{R} (\varepsilon) h^{\beta})} |
  \nabla V (x) - \nabla V ( \hat{q}_h (t)) |u_h^2 (t, x)
  dx}{\int_{\mathbb{R}^N} u_h^2 (t, x) dx} ;\\
  L_2 & = & \frac{\int_{\mathbb{R}^N \smallsetminus B ( \hat{q}_h (t), \hat{R}
  (\varepsilon) h^{\beta})} | \nabla V (x) |u_h^2 (t, x)
  dx}{\int_{\mathbb{R}^N} u_h^2 (t, x) dx} ;\\
  L_3 & = & \frac{\int_{\mathbb{R}^N \smallsetminus B ( \hat{q}_h (t), \hat{R}
  (\varepsilon) h^{\beta})} | \nabla V ( \hat{q}_h (t)) |u_h^2 (t, x)
  dx}{\int_{\mathbb{R}^N} u_h^2 (t, x) dx} .\\
  &  &
\end{eqnarray*}
By the Theorem \ref{teoconcinf} and by Lemma \ref{lemmabar5} we have $L_3 < M
\varepsilon$.

We have also
\begin{equation}
L_1 \leq K_1 + R_2 + 1.
\end{equation}
Using hypothesis (\ref{Vinf}) we have
\begin{equation}
  L_2 \leq M \varepsilon + \left[ \frac{L}{\sigma^2} \right]^b \varepsilon^{1- b},
\end{equation}
where $b \in \left( 0, 1 \right)$ is defined in (\ref{Vinf}).
Concluding we have
\begin{equation}
  |H_h (t) | \leq \frac{3 LM}{\sigma^2 R_3^{a - 1}} + \left[
  \frac{L}{\sigma^2} \right]^b \varepsilon^{1 - b} + M (2 + 3 R_3) \varepsilon
  + 2 M \hat{R} (\varepsilon) h^{\beta} .
\end{equation}
At this point we can have $\sup_t |H_h (t) |$ arbitrarily small choosing
firstly $R_3$ sufficiently large, secondly $\varepsilon$ sufficiently small,
and finally $h$ small enough.$\square$
\begin{corollary}
Let $\psi (t, x) \in C (\mathbb{R}, H^2) \cap C^1 ( \mathbb{R}, H^1)$
  be a global solution of
  equation (\ref{schr}), with initial data
 $\varphi(x)=U(\frac{x-q_0}{h^\beta})e^{\frac ih \mathbf{v} \cdot x}$ where
  $U$ is a radial minimizer of $J$ on $S\sigma$, 
 $q_0\in\R^N$, $\mathbf{v}\in\R^N$, and
  $h<h_0$. Then the barycenter $q$ satisfies the following Cauchy problem
\begin{equation*}
\left\{
\begin{array}{l}
\ddot{q}_{h}(t)+\nabla V(q_{h}(t))=H_{h}(t) \\
q_{h}(0)=q_{0} \\
\dot{q}_{h}(0)=\mathbf{v}%
\end{array}%
\right.
\end{equation*}
\end{corollary}
\noindent\textbf{Proof.} The initial data belongs to $B_h^K$ for some $K$.
We apply the previous results to
obtain the equation for $\ddot q$. The initial data $q(0)$ and $\dot q(0)$ are derived with a
direct calculation. $\square$

\section{The swarm interpretation}

In this section we present a different point of view on our problem.
Although this approach is non rigorous, it provides some physical intuitions which are
inspiring for a better understanding of the general framework.
We will suppose that the soliton is composed by a swarm of
particles which follow the laws of classical dynamics given by the
Hamilton-Jacobi equation. This interpretation will permit us to give an
heuristic proof of the main result.

First of all let us write NSE with the usual physical constants $m$ and
$\hslash$:
\[ i \hslash \frac{\partial \psi}{\partial t} = - \frac{\hslash^2}{2 m} \Delta
   \psi + \frac{1}{2} W_{\hslash}' (\psi) + V (x) \psi . \]
Here $m$ has the dimension of {\em{mass}} and $\hslash$, the Plank
constant, has the dimension of {\em{action}}.

In this case equations (\ref{Sh1}) and (\ref{Sh2}) become:
\begin{equation}
  - \frac{\hslash^2}{2 m} \Delta u + \frac{1}{2} W_{\hslash}' (u) + \left(
  \partial_t S + \frac{1}{2 m} \left| \nabla S \right|^2 + V (x) \right) u = 0
  ; \label{Sh1c}
\end{equation}
\begin{equation}
  \partial_t \left( u^2 \right) + \nabla \cdot \left( u^2 \frac{\nabla S}{m}
  \right) = 0. \label{Sh2c}
\end{equation}
The second equation allows us to interpret the matter field to be a fluid
composed by particles whose density is given by
\[ \rho_{\mathcal{H}} = u^2 \]
and which move in the velocity field
\begin{equation}
  \mathbf{v} = \frac{\nabla S}{m} . \label{linda}
\end{equation}
So equation (\ref{Sh2c}) becomes the continuity equation:
\[ \partial_t \rho_{\mathcal{H}} + \nabla \cdot \left( \rho_{\mathcal{H}}
   \mathbf{v} \right) = 0. \]
If
\begin{equation}
  - \frac{\hslash^2}{2 m} \Delta u + \frac{1}{2} W_{\hslash}' (u) \ll u,
  \label{rosaS}
\end{equation}
equation (\ref{Sh1c}) can be approximated by the eikonal equation
\begin{equation}
  \partial_t S + \frac{1}{2 m} \left| \nabla S \right|^2 + V (x) = 0.
  \label{hjS}
\end{equation}
This is the Hamilton-Jacobi equation of a particle of mass $m$ in a potential
field $V$.

If we do not assume (\ref{rosaS}), equation (\ref{hjS}) needs to be replaced
by
\begin{equation}
  \partial_t S + \frac{1}{2 m} \left| \nabla S \right|^2 + V + Q (u) = 0
  \label{hjqS}
\end{equation}
with
\[ Q (u) = \frac{- \left( \hslash^2 / m \right) \Delta u + W_{\hslash}' (u)}{2
   u} . \]
The term $Q (u)$ can be regarded as a field describing a sort of interaction
between particles.

Given a solution $S (t, x)$ of the Hamilton-Jacobi equation, the motion of the
particles is determined by Eq.(\ref{linda}).

\subsection{An heuristic proof}

In this section we present an heuristic proof of the main result.
This proof is not at all rigorous, but it helps to understand the
underlying Physics.

\

If we interpret $\rho _{\mathcal{H}}=u^{2}$ as the density of
particles then
\begin{equation*}
\mathcal{H=}\int \rho _{\mathcal{H}}dx
\end{equation*}%
is the total number of particles. By (\ref{hjqS}), each of these particle
moves as a classical particle of mass $m$ and hence, we can apply to the
laws of classical dynamics. In particular the center of mass defined in
(\ref{bary}) takes the following form:
\begin{equation}
q(t)=\frac{\int xm\rho _{\mathcal{H}}dx}{\int m\rho _{\mathcal{H}}dx}=
\frac{\int x\rho _{\mathcal{H}}dx}{\int \rho _{\mathcal{H}}dx}.  \label{aiace}
\end{equation}
The motion of the barycenter is not affected by the interaction between
particles (namely by the term (\ref{hjqS})), but only by the external
forces, namely by $\nabla V.$ Thus the global external force acting on the
swarm of particles is given by
\begin{equation}
\overrightarrow{F}=-\int \nabla V(x)\rho _{\mathcal{H}}dx.  \label{ulisse}
\end{equation}
Thus the motion of the center of mass $q$ follows the Newton law
\begin{equation}
\overrightarrow{F}=M{\ddot{q}},  \label{tersite}
\end{equation}
where $M=\int m\rho _{\mathcal{H}}dx$ is the total mass of the swarm; thus
by (\ref{aiace}), (\ref{ulisse}) and (\ref{tersite}), we get
\begin{equation*}
{\ddot{q}}(t)=-\frac{\int \nabla V\rho _{\mathcal{H}}dx}
{m\int \rho _{\mathcal{H}}dx}=-\frac{\int \nabla Vu^{2}dx}{m\int u^{2}dx}.
\end{equation*}

If we assume that the $u(t,x)$ and hence $\rho _{\mathcal{H}}(t,x)$ is
concentrated in the point $q(t),$ we have that
\begin{equation*}
\int \nabla Vu^{2}dx\cong \nabla V\left( q(t)\right) \int u^{2}dx
\end{equation*}
and so, we get
\begin{equation*}
m{\ddot{q}}(t)\cong -\nabla V\left( q(t)\right).
\end{equation*}

Notice that the equation $m{\ddot{q}}(t)=-\nabla V\left( q(t)\right) $ is
the Newtonian form of the Hamilton-Jacobi equation (\ref{hjS}).

\end{document}